\newif\ifSHOWEXTRA
\SHOWEXTRAtrue
\SHOWEXTRAfalse
\documentclass[12pt]{article}
\usepackage[utf8]{inputenc}
\usepackage[margin=1in]{geometry}
\usepackage{amsmath,amsthm,amssymb,amsfonts, pifont, bbm, xcolor, scrextend}
\definecolor{forestgreen}{RGB}{0,98,51}
\usepackage{appendix}
\usepackage{fancyhdr}
\usepackage{color,colortbl}
\newcounter{CommentCounter}
\newcommand{\daniel}[1]{{\color{purple} \sf $\spadesuit\spadesuit\spadesuit$ Daniel (\stepcounter{CommentCounter}\theCommentCounter): [#1]}}

\definecolor{Gray}{gray}{0.9}
\definecolor{LightGreen}{rgb}{0,0.1,0.1}
\newcommand{\R}{\mathbb R}

\newcommand{\Pa}{\mathbb P}
\newcommand{\Ea}{\mathbb E}

\newcommand{\N}{\mathbb N}
\newcommand{\Q}{\mathbb Q}
\newcommand{\Z}{\mathbb Z}

\usepackage{tikz}
\usepackage{pgfplots}
\usepackage{eqnarray}
\usepackage{outlines}
\usepackage[mathscr]{euscript}
 \let\mathscr\relax% just so we can load this and rsfs
\usepackage[scr]{rsfso}
\usepackage{multicol}
\usepackage{enumerate}
\usepackage[colorlinks=true, pdfstartview=FitV, linkcolor=blue,
citecolor=blue, urlcolor=blue]{hyperref}

\newcommand{\BIX}{\overline{\bf X}}
\newcommand{\biX}{\overline{X}}

\newtheorem{thm}{Theorem}[section]
\numberwithin{thm}{section}

\newtheorem{prop}[thm]{Proposition}
\newtheorem{lem}[thm]{Lemma}

\numberwithin{quest}{section}

\theoremstyle{remark}
\newtheorem{rem}{Remark}
\numberwithin{rem}{section}

\theoremstyle{definition}

\numberwithin{claim}{thm}

\pgfplotsset{compat=1.12}

\begin{document}
\author{
Daniel J. Slonim\footnote{Purdue University,
  Department of Mathematics,
150 N. University Street, West Lafayette, IN 47907, dslonim@purdue.edu, 0000-0002-9554-155X}
}
\title{Ballisticity of Random walks in Random Environments on $\Z$ with Bounded Jumps}

\date{\today}

\maketitle

\begin{abstract}
    We characterize ballistic behavior for general i.i.d. random walks in random environments on $\Z$ with bounded jumps. The two characterizations we provide do not use uniform ellipticity conditions. They are natural in the sense that they both relate to formulas for the limiting speed in the nearest-neighbor case. 
    \medskip\noindent {\it MSC 2020itions. :}
    60G50 %Sums of independent random variables; random walks
    60J10 %Markov Chains  (discrete-time Markov processes on discrete state spaces)
    60K37 %Processes in Random Environments
    \\
    {\it Kewords:}
    random walk,
    random environment,
    bounded jumps,
    ballisticity
\end{abstract}

%\tableofcontents

\section{Introduction}

In this paper, we provide two characterizations of ballisticity for random walks in random environments (RWRE) on $\Z$ with bounded jumps. Most previous characterizations of ballisticity for such RWRE (or for RWRE on a strip, a generalization of the one-dimensional bounded-jump model) are in terms of limits of norms of products of random matrices that are difficult or impossible to check in practice, and involve strong ellipticity assumptions that preclude certain types of environments. (See, for example, \cite{Bremont2002}, \cite{Bremont2004}, \cite{Goldsheid2008}, \cite{Roitershtein2008}, \cite{Bremont2009}). 

Our characterizations of ballisticity are intuitively very easy to understand (if not to check in general), and do not use strong ellipticity assumptions. The primary motivation for our characterizations is that although we do not have a way to check them in general, we are able to check them in the case of Dirichlet environments, a special, weakly elliptic model of RWRE where certain exact computations are often possible (sometimes not for the walk one wants to study, but for walks on finite graphs that can be related to the desired model). The author does this in \cite{Slonim2021a}. Here, we are focused only on the general case. 

\subsection{Model}\label{subsec:model}
An {\em environment on $\Z$} is a nonnegative function $\omega:\Z\times \Z\to [0,1]$ such that for all $x\in \Z$, $\sum_{y\in \Z}\omega(x,y)=1$. For a fixed $x$ and $\omega$, we will let $\omega^x$ be the measure on $\Z$ defined by $\omega^x(y)=\omega(x,x+y)$. Then we can identify the function $\omega$ with the tuple $(\omega^x)_{x\in\Z}$. Let $\mathcal{M}_1(\Z)$ be the set of probability measures on $\Z$,\ (endowed with the topology of weak convergence); then $\Omega:=\prod_{x\in \Z}\mathcal{M}_1(\Z)$ is the set of all environments on $\Z$. 

For a given environment $\omega$ and $x\in \Z$, we can define the {\em quenched} measure $P_{\omega}^x$ on $\Z^{\N}$ (where we assume $0\in \N$) to be the law of a Markov chain ${\bf X}=(X_n)_{n\geq0}$ on $\Z$, started at $x$, with transition probabilities given by $\omega$. That is, $P_{\omega}^x(X_0=x)=1$, and for $n\geq1$, $P_{\omega}^x(X_{n+1}=y|X_0,\ldots,X_n)=\omega(X_n,y)$. 

Let $\mathcal{F}$ be the Borel sigma field with respect to the product topology on $\Omega$,
and let $P$ be a probability measure on $(\Omega,\mathcal{F})$. For a given $x\in \Z$, we define the {\em annealed} measure $\Pa^x=P\times P_{\omega}^x$ on $\Omega\times \Z^{\N}$ by
\begin{equation*}\Pa^x(A\times B)=\int_{A}P_{\omega}^x(B)P(d\omega)\end{equation*}
for measurable $A\subset\Omega,B\subset \Z^{\N}$. In particular, for all measurable events $B\subset \Z^{\N}$, we have $\Pa^x(\Omega\times B)=E[P_{\omega}^x(B)]$. We often abuse notation by writing $\Pa^x(B)$ instead of $\Pa^x(\Omega\times B)$. 

As another notational convenience, we will use interval notation to denote sets of consecutive integers in the state space $\Z$, rather than subsets of $\R$. Thus, for example, we will use $[1,\infty)$ or $(0,\infty)$ to denote the set of integers strictly to the right of 0. However, we make one exception, using $[0,1]$ to denote the set of all real numbers from 0 to 1.  

For a subset $S\subseteq\Z$, let $\omega^S=(\omega^x)_{x\in S}$. In the case where $S$ is a half-infinite interval, we simplify our notation by using $\omega^{\leq x}$ to denote $\omega^{(-\infty,x]}$, and similarly with $\omega^{<x}$, $\omega^{\geq x}$, and $\omega^{>x}$. 

We consider the following conditions for a probability measure $P$ on $\Omega$:
\begin{enumerate}[(C1)]
    \item\hypertarget{cond:C1} The $\{\omega^x\}_{x\in\Z}$ are i.i.d. under $P$.
    \item\hypertarget{cond:C2} For $P$--a.e. environment $\omega$, the Markov chain induced by $P_{\omega}^0$ is irreducible.
    \item\hypertarget{cond:C3} There exist $L$ and $R$ such that for $P$--a.e. environment $\omega$, $\omega(a,b)=0$ whenever $b$ is outside $[a-L,a+R]$. 
\end{enumerate}

\subsection{Results}
It was shown in \cite{Key1984} that under the above assumptions, a 0-1 law holds for directional transience. That is, the walk is either almost surely transient to the right, almost surely transient to the left, or almost surely recurrent. \ifSHOWEXTRA
\daniel{His middle assumption wasn't irreducibility, it was just that at every site there's a.s. a positive probability of stepping to the left and the right.}
\fi
We can also show that under these assumptions, a limiting velocity necessarily exists. 
\begin{prop}\label{prop:LimitingVelocity} Let $P$ be a probability measure on $\Omega$ satisfying (\hyperlink{cond:C1}{C1}), (\hyperlink{cond:C2}{C2}), and (\hyperlink{cond:C3}{C3}). Then there is a $\Pa^0$--almost sure limiting velocity $v=\lim_{n\to\infty}\frac{X_n}{n}$. Moreover, $\lim_{x\to\infty}\frac{H_{\geq x}}{x} = \frac{1}{v}$, where $\frac{1}{v}$ is understood to be $\infty$ if $v=0$.
\end{prop}
It was seen in \cite{Bremont2004} that this limiting velocity exists under a uniform ellipticity assumption, 
\ifSHOWEXTRA
\daniel{\cite{Bremont2004} uses uniform ellipticity, but \cite{Bremont2002}, which assumes $R=1$, uses a middle-ground assumption. The paper \cite{Bremont2004} simply cites the argument in \cite{Bremont2002}, which in turn pretty much just appeals to Kozlov's arguments on invariant measure as best as I can tell. It's not easy to tell where, if at all, the ellipticity assumptions come in without thoroughly reading multiple papers.} 
\fi
but it can be proven in the more general case with standard techniques, which we outline in Section \ref{subsec:existence}.
We then provide a characterization of ballisticity, making the following additional assumption for convenience.
\begin{enumerate}[(C4)]
    \item\hypertarget{cond:C4} For $P$--a.e. environment $\omega$, $\lim_{n\to\infty}X_n=\infty$, $P_{\omega}^0$--a.s.
\end{enumerate}
By symmetry, our characterization also handles the case where the walk is transient to the left, and thus by the 0-1 law for directional transience, completely characterizes the regime $v\neq0$ for all measures $P$ satisfying (\hyperlink{cond:C1}{C1}), (\hyperlink{cond:C2}{C2}), and (\hyperlink{cond:C3}{C3}).

To formally state our characterization, we must establish notation for hitting times, as well as notation that counts the number of visits to a given site. For a given walk ${\bf X}=(X_n)_{n=0}^{\infty}$, we define $H_x({\bf X})$ to be the first time the walk hits $x\in\Z$. That is,
\begin{equation*}
    H_x({\bf X})=\inf\{n\in\N:X_n=x\}
\end{equation*}
We usually write it as $H_x$ when we can do so without ambiguity. For a subset $S\subset\Z$, let $H_S=\min_{x\in S}H_x$. First positive hitting times are denoted as $\tilde{H}_x$ or $\tilde{H}_S$. That is,
\begin{equation*}
    \tilde{H}_x({\bf X})=\inf\{n\in\N:X_n=x\},
\end{equation*}
and $\tilde{H}_S=\min_{x\in S}\tilde{H}_x$.
If the set is the half-infinite interval $[x,\infty)$, we use $H_{\geq x}$ to denote its hitting time, and similarly with $H_{>x}$, $H_{\leq x}$, and $H_{>x}$. 
For a walk ${\bf X} = (X_n)_{n=0}^{\infty}$ on $\Z$ with $x\in \Z$, $N_x({\bf X})=\#\{n\in\N:X_n=x\}$ is the number of times the walk is at site $x$. We usually write it as $N_x$ if we are able to do so without ambiguity. For a subset $S\subset\Z$, let $N_S=\sum_{x\in S}N_x$. 

\begin{thm}\label{thm:abstractballisticitycriteria}
Let $P$ be a probability measure on $\Omega$ satisfying (\hyperlink{cond:C1}{C1}), (\hyperlink{cond:C2}{C2}), (\hyperlink{cond:C3}{C3}), and (\hyperlink{cond:C4}{C4}). Then the following are equivalent:
\begin{enumerate}[(a)]
\item The walk is ballistic: $v>0$.
\item $\Ea^0[H_{\geq1}]<\infty$.
\item $\Ea^{0}[N_{0}]=E[E_{\omega}^0[N_0]]<\infty$.
\end{enumerate} 
\end{thm}

\begin{rem}
The equivalence of (a) and (b) was proven under a strong ellipticity assumption by Brémont (see \cite[Theorem 3.7]{Bremont2002}, \cite[Proposition 9.1]{Bremont2004}). We prove the equivalence of (a), (b), and (c) without such an assumption. 
\end{rem}

\begin{rem}
These characterizations are quite natural, given that in the nearest-neighbor case we in fact have the identity
\begin{equation}
v=\frac{1}{\Ea^0[N_0]}=\frac{1}{\Ea^0[H_{\geq1}]}.
\end{equation}
In general, the above turns into
\begin{equation}
v\geq\frac{1}{\Ea^0[N_0]}\geq\frac{1}{\Ea^0[H_{\geq1}]}.
\end{equation}
as we shall see from Lemmas \ref{prop:4.8}, \ref{lem:LimitingVelocity}, and \ref{lem:limitfinitecase}.
\end{rem}

%\begin{rem}
%An analog of Theorem \ref{thm:abstractballisticitycriteria} can be proven for RWRE on a strip (provided certain conditions are satisfied), but we confine our discussion to RWRE on $\Z$ with bounded jumps for simplicity and because no real new insight seems to come from doing this work in the strip model. \daniel{Look into this.}
%\end{rem}

\section{Proofs}

%To do this, we first discuss a way to generate a ``walk from $-\infty$ to $\infty$" in a given environment with transience to the right. This walk exhibits a sort of stationary behavior, and we show that the limiting velocity of a walk started at 0 is the reciprocal of the annealed expected amount of time this bi-infinite walk spends at a given site. While we do not know how to calculate this quantity, we show that it is finite if the annealed expected amount of time walk started at 0 spends at 0 is finite. 

%\begin{prop}\label{prop:positivespeed}
%Let $P$ be a probability measure on $\Omega$ satisfying (\hyperlink{cond:C1}{C1}), (\hyperlink{cond:C2}{C2}), (\hyperlink{cond:C3}{C3}), and (\hyperlink{cond:C4}{C4}). Let $\accel$ be a Bouchet acceleration function. Then If $\Ea^0[N_0]<\infty$, then $v>0$.
%\end{prop}

%\ifSHOWEXTRA
%\hyperlink{proof:positivespeed}{Jump to proof.}
%\fi

%Our results characterizing ballisticity require only the treatment of the discrete-time case. However, in Section \ref{sec:acceleration}, we will analyze the ballisticity of an accelerated walk in continuous time, following \cite{Bouchet2013}. For that reason, we want to treat the continuous case as well. The discrete-time case may be viewed as a special case of the continuous case in which the jump times are all 1 almost surely. Thus, we assume the continuous case throughout this subsection.
%However, because the discrete-time case is the main focus of our attention, we will continue to use $n$ rather than $t$ as the index of a walk where applicable.

We discuss the proof of Proposition \ref{prop:LimitingVelocity}, and then provide a full proof of Theorem \ref{thm:abstractballisticitycriteria}. 

\subsection{Existence of limiting velocity}\label{subsec:existence}

Becasue the proof of Proposition \ref{prop:LimitingVelocity} is only a slight modification of work that has already been done, we simply outline some details of the argument rather than giving a full proof. 

The proof for the recurrent case (where, necessarily, $v=0$) can be done by a slight modification of arguments in \cite{Zerner2002}, which we leave to the reader. The proof for the directionally transient case follows \cite{Kesten1977} in defining {\em regeneration times} $(\tau_k)_{k=0}^{\infty}$. Let $\tau_0:=0$, and for $k\geq1$, define 
\begin{equation}\label{eqn:RegenerationTimes}
    \tau_k:=\min\{n>\tau_{k-1}:X_n>X_j\text{ for all }j<n,X_n\leq X_j\text{ for all }j>n\}.
\end{equation}
We can then show
\begin{equation}\label{eqn:856}
    v:=\lim_{n\to\infty}\frac{X_n}{n}=\frac{\Ea^0[X_{\tau_{2}}-X_{\tau_1}]}{\Ea[\tau_{2}-\tau_{1}]},
\end{equation}
where the numerator is always finite and the fraction is understood to be 0 if the denominator is infinite.
It is standard (e.g., \cite{Kesten1977}, \cite{Sznitman&Zerner1999}) to prove a LLN like \eqref{eqn:856} by the following steps:
\begin{enumerate}[(a)]
    \item Show that $\frac{X_{\tau_k}}{k}$ approaches $\Ea[X_{\tau_2}-X_{\tau_1}]$.
    \item Show that $\frac{\tau_k}{k}$ approaches $\Ea[\tau_2-\tau_1]$.
    \item Show that $\Ea[X_{\tau_2}-X_{\tau_1}]<\infty$.
    \item Conclude that the limit \eqref{eqn:856} holds for the subsequence $\left(\frac{X_{\tau_k}}{\tau_k}\right)_k$.
    \item Use straightforward bounds that come from the definitions of the $\tau_k$ to get the limit for the entire sequence $\left(\frac{X_n}{n}\right)_n$. 
\end{enumerate}
The identity $\lim_{x\to\infty}\frac{H_{\geq x}}{x} = \frac{1}{v}$ then comes from a comparison of $\frac{x}{H_{\geq x}}$ with a subsequence of $\frac{X_n}{n}$. 

The definition of the regeneration times is precisely set up so that both the sequences $(\tau_k-\tau_{k-1})_{k\geq2}$ and $(X_{\tau_k}-X_{\tau_{k-1}})_{k\geq2}$ are i.i.d. sequences, so proving the limits (a) and (b) is a matter of tracing how the i.i.d. feature follows from the definitions and applying the strong law of large numbers. In fact, arguing as in \cite[Lemma 1]{Kesten1977}, one can show that the triples
\begin{equation*}
    \xi_n:=\left(\tau_n-\tau_{n-1}~,~(X_{\tau_{n-1}+i}-X_{\tau_{n-1}})_{i=1}^{\tau_n-\tau_{n-1}}~,~(\omega^x)_{x=X_{\tau_{n-1}}}^{X_{\tau_n}-1}\right)
\end{equation*}
are i.i.d. under $\Pa^0=P\times P_{\omega}^0$ for $n\geq2$.

The finiteness in (c) can be shown using arguments along the lines of those in \cite[Lemma 3.2.5]{Zeitouni2004}. Because we have assumed almost-sure transience to the right, the measure $\Q$ introduced there is unnecessary. Another difference is that in our model, transience to the right does not imply that every vertex to the right of the origin is hit. There is a point in the argument from \cite{Zeitouni2004} where Zeitouni argues the the $\Q$-probability, for a given $x$, that a regeneration occurs at site $x$ approaches $\Pa^0(H_{<0}=\infty)$. Instead, we focus on the probability that the regeneration occurs on a given interval of length $R$. For $z\geq0$, let $B_z$ be the event that for some $k$, $X_{\tau_k}\in[zR,(z+1)R)$. Then
\begin{align}
    \notag
    \Pa^0(B_z)&=E\left[P_{\omega}^0(B_z)\right]
    \\
    \notag
    &\geq
    E\left[\sum_{i=0}^{R-1}P_{\omega}^0(X_{H_{[zR,(z+1)R)}}=zR+i)P_{\omega}^{zR+i}(H_{<zR+i}=\infty)\right]
    \\
    \notag
    &= 
    \sum_{i=0}^{R-1}E\left[P_{\omega}^0(X_{H_{[zR,(z+1)R)}}=zR+i)P_{\omega}^{zR+i}(H_{<zR+i}=\infty)\right]
    \\
\notag
    &=\sum_{i=0}^{R-1}\Pa^0(X_{H_{[zR,(z+1)R)}}=zR+i)\Pa^{zR+i}(H_{<zR+i}=\infty)
    \\
    \label{eqn:1935}
    &=\Pa^{0}(H_{<0}=\infty),
\end{align}
where the second to last equality comes from the fact that $\omega^{<zR}$ is independent of $\omega^{\geq zR+i}$, and the last comes from translation invariance and the fact that $H_{[zR,(z+1)R)}<\infty$ $\Pa^0$--a.s. The rest of the argument from \cite{Zeitouni2004} goes through to prove (c), and (d) and (e) easily follow.

\subsection{Ballisticity}

For the rest of this paper, assume $P$ satisfies (\hyperlink{cond:C1}{C1}), (\hyperlink{cond:C2}{C2}), (\hyperlink{cond:C3}{C3}), and (\hyperlink{cond:C4}{C4}). Our goal is to prove Theorem \ref{thm:abstractballisticitycriteria}.
We begin with the following lemma, from which $(b)\Rightarrow(c)$ follows immediately.

\begin{lem}\label{prop:4.8}
$\Ea^0[N_0]\leq\Ea^0[H_{\geq1}]$. 
\end{lem}
\begin{proof}
The visits to 0 may be sorted based on the farthest point to the right that the walk has hit in the past at the time of each visit. For a given $y<x\in\Z$, define $N_{y}^{(-\infty,x)}$ to be the amount of time the walk spends at $y$ before $H_{\geq x}$..
Thus, for a walk started at 0 we get
\begin{equation}\label{eqn:1207}
    N_0=\sum_{x=0}^{\infty}\left(N_0^{(-\infty,x+1)}-N_0^{(-\infty,x)}\right).
\end{equation}
Taking expectations on both sides, we get
\begin{equation}\label{eqn:1213}
    \Ea^0[N_0]
    =\sum_{x=0}^{\infty}E\left[E_{\omega}^0\left[N_0^{(-\infty,x+1)}-N_0^{(-\infty,x)}\right]\right]
\end{equation}
Now $N_0^{(-\infty,x)}$ and $N_0^{(-\infty,x+1)}$ can only differ if the walk hits $[x,\infty)$ at $x$. Conditioned on this event, the distribution under $P_{\omega}^0$ of %the walk
$(X_{n+H_{\geq x}})_{n=0}^{\infty}$ is the distribution of ${\bf X}$ under $P_{\omega}^x$. Thus, 
\begin{equation}\label{eqn:1218}
    E_{\omega}^0\left[N_0^{(-\infty,x+1)}-N_0^{(-\infty,x)}\right]
    =
    P_{\omega}^0(X_{H_{\geq x}}=x)E_{\omega}^x\left[N_0^{(-\infty,x+1)}\right].
\end{equation}
Combining \eqref{eqn:1213} and \eqref{eqn:1218}, we get
\begin{align*}
    \Ea^0[N_0]
    &=\sum_{x=0}^{\infty}E\left[P_{\omega}^0(X_{H_{\geq x}}=x)E_{\omega}^x\left[N_0^{(-\infty,x+1)}\right]\right]
    \\
    &\leq \sum_{x=0}^{\infty}E\left[E_{\omega}^x\left[N_0^{(-\infty,x+1)}\right]\right]
    \\
    &=\sum_{x=0}^{\infty}\Ea^x\left[N_0^{(-\infty,x+1)}\right].
\end{align*}
By stationarity, 
\begin{align*}
    \Ea^0[N_0]
    &\leq\sum_{x=0}^{\infty}\Ea^0\left[N_{-x}^{(-\infty,1)}\right]
    \\
    &=\Ea^0\left[\sum_{x=0}^{\infty}N_{-x}^{(-\infty,1)}\right]
    \\
    &=\Ea^0[H_{\geq1}].
\end{align*}
This completes the proof.
\end{proof}

Our next goal is to prove (c) $\Rightarrow$ (a). 
The proof of the following lemma, will use the regeneration times defined in \eqref{eqn:RegenerationTimes}.

\begin{lem}\label{known}
For any $c\in\Z$,
\begin{equation*}\lim_{x\to\infty}\frac{1}{x}\sum_{y=c}^xN_y=\frac1{v},~\Pa^a\text{--a.s.}\end{equation*}
If $v=0$, then the limit is infinity.
\end{lem}

\begin{proof}
Fix $c$. 
As in the proof of Lemma \ref{prop:4.8}, we use $N_{y}^{(-\infty,x)}$ to denote the amount of time the walk spends at $y$ before $H_{\geq x}$.
Then for any $x>c$, write

%\begin{equation*}H_{\geq x}=\sum_{y=-\infty}^{c-1}N_y^{(-\infty,x)}+\sum_{y=c}^{x-1}N_y^{(-\infty,x)}\end{equation*}

\begin{equation*}\frac{H_{\geq x}}{x}=\frac1x\sum_{y=-\infty}^{c-1}N_y^{(-\infty,x)}+\frac1x\sum_{y=c}^{x-1}N_y^{(-\infty,x)}.
\end{equation*}
The first sum is bounded above by $\sum_{y=-\infty}^{c-1}N_y$, which is almost surely finite by assumption (\hyperlink{cond:C4}{C4}). Dividing by $x$ therefore sends the first term to 0 as $x\to\infty$; hence, by Proposition \ref{prop:LimitingVelocity},

\begin{equation}\label{eqn:890}
\lim_{x\to\infty}\frac1x\sum_{y=c}^{x-1}N_y^{(-\infty,x)}=\frac1{v}, \Pa^a\text{--a.s.}
\end{equation}

We note that $N_y$ and $N_y^{(-\infty,x)}$ differ only if the walk backtracks and visits $y$ after reaching $[x,\infty)$. The sum, over all $y<x$, of these differences, is the total amount of time the walk spends to the left of $x$ after $H_{\geq x}$, and it is bounded above by the time from $H_{\geq x}$ to the next regeneration time (defined as in \eqref{eqn:RegenerationTimes}), which is in turn bounded above by $\tau_{J(x)}-\tau_{J(x)-1}$, where $J(x)$ is the (random) $j$ such that $\tau_{j-1}\leq H_{\geq x}<\tau_j$. Hence

\begin{equation}\label{eqn:896}
\frac1x\sum_{y=c}^{x-1}N_y^{(-\infty,x)}
\leq \frac1x\sum_{y=c}^{x-1}N_y
\leq \frac1x\sum_{y=c}^{x-1}N_y^{(-\infty,x)} + \frac1x[\tau_{J(x)}-\tau_{J(x)-1}].
\end{equation}

Assume $v=0$. Then by \eqref{eqn:890}, the left side of \eqref{eqn:896} approaches $\infty$ as $x$ approaches $\infty$, and therefore so does the middle. 
On the other hand, suppose $v>0$. By \eqref{eqn:856}, $\Ea[\tau_{2}-\tau_{1}]<\infty$.
Then by the strong law of large numbers, $\frac{\tau_n}{n}\to \Ea[\tau_{2}-\tau_{1}]<\infty$, which implies that $\frac{\tau_n-\tau_{n-1}}{n}$ approaches 0. Since $J(x)\leq x+1$, the term $\frac1x[\tau_{J(x)}-\tau_{J(x)-1}]$ approaches zero almost surely; hence the squeeze theorem yields the desired result.
\end{proof}

Next, we will define a ``walk from $-\infty$ to $\infty$."
%Suppose for now that $R=1$. Then, for almost every $\omega$, it is easy to define a bi-infinite walk $\BIX=(\biX_n)_{n\in\Z}$ whose ``right halves" are distributed like random walks under $\omega$. From each site $a$, run a walk according to the transition probabilities given by $\omega$ until it reaches $a+1$ (which occurs in finite time $P_{\omega}^a$--a.s. for $P$--a.e. $\omega$). Concatenating all of these walks then gives, up to a time shift\footnote{Choose, for example, the time shift where $\biX_0=0$ and where $\biX_n<0$ whenever $n<0$.}, a unique walk $\BIX=(\biX_n)_{n\in\Z}$ such that for any $x\in\Z$, the distribution of $(\biX_k)_{k=n}^{\infty}$, conditioned on $\biX_n=x$, is $P_{\omega}^x$. We may think of $\BIX$ as a walk from $-\infty$ to $\infty$ in the environment $\omega$.
%
%With a bit more work, we can define a similar bi-infinite walk in the general case $R>0$. 
Call the set of vertices $((k-1)R,kR]$ the $k$th {\em level} of $\Z$, and for $x\in\Z$, let $[[x]]_R$ denote the level containing $x$. Let $\omega$ be a given environment. From each point $a\in\Z$, run a walk according to the transition probabilities given by $\omega$ until it reaches the next level (i.e., $[[a+R]]_R$). This will happen $P_{\omega}^a$--a.s. for $P$--a.e. $\omega$, by assumption (\hyperlink{cond:C4}{C4}) and because it is not possible to jump over a set of length $R$. Do this independently at every point for every level. This gives what we will call a {\em cascade}: a set of (almost surely finite) walks indexed by $\Z$, where the walk indexed by $a\in\Z$ starts at $a$ and ends upon reaching level $[[a+R]]_R$. Equip the set of cascades with the natural sigma field, let $P_{\omega}$ be the probability measure we have just described on the space of cascades, and let $\Pa=P\times P_{\omega}$. 

For $\Pa$--almost every cascade (i.e., those where the walk started from each vertex hits the level to its right), we can concatenate an appropriate chain of these finite walks to generate a walk started at any point $a\in\Z$. To the walk started from $a=a_0$, append the walk started from the point $a_1\in[[a+R]]_R$ where that walk lands. And to that, append the walk started from the point $a_2\in[[a+2R]]_R$ where the finite walk from $a_1$ lands, and so on. This gives, for each point $a$, a right-infinite walk ${\bf X}^{a}=(X_n^a)_{n=0}^{\infty}$. 
It is crucial to note that by the strong Markov property, the law of ${\bf X}^{a}$ under $P_{\omega}$ is the same as the law of ${\bf X}$ under $P_{\omega}^{a}$, which also implies that the law of ${\bf X}^{a}$ under $\Pa$ is the same as the law of ${\bf X}$ under $\Pa^{a}$.

%Now if there is only one path, $N_{x}$ denotes the number of times that one path ${\bf X}$ hits the point $x$. Hence $N_x$ is a function defined on the space of right-infinite paths. On the other hand, we may define a function $N_{x}^{a}$ on the space of cascades to be the number of times the path ${\bf X}^{a}$ hits the point $(x,y)$. Again, we may note, for example, that ${\Ea\left[N_x^{a}\right]=\Ea^{a}\left[N_x\right]}$.

For each $x\in\Z$, let the ``coalescence event" $C_x$ be the event that all the walks from level $[[x-R]]_R$ first hit level $[[x]]_R$ at $x$. On the event $C_x$, we say a coalescence occurs at $x$.

\begin{lem}\label{prop:coalescences}
Let $\mathcal{E}_1$ be the event that all the ${\bf X}^a$ are transient to the right, that all steps to the left and right are bounded by $L$ and $R$, respectively, and that infinitely many coalescences occur to the left and to the right of 0. Then $\Pa(\mathcal{E}_1)=1$.
\end{lem}

\begin{proof}
Boundedness of steps has probability 1 by assumption (\hyperlink{cond:C3}{C3}), and by assumption (\hyperlink{cond:C4}{C4}) all the walks ${\bf X}^a$ are transient to the right with probability 1. Now for $k\geq 2$ and $x\in\Z$, let $C_{x,k}$ be the event that all the walks from level $[[x-R]]_R$ first hit level $[[x]]_R$ at $x$ without ever having reached $[[x-kR]]_R$. Choose $k$ large enough that $\Pa(C_{0,k})>0$; then under the law $\Pa$, the events $\{C_{nkR,k}\}_{n\in\Z}$ are all independent and have equal, positive probability. Thus, infinitely many of them will occur in both directions, $\Pa$--a.s. By definition, $C_{x,k}\subset C_x$, and so infinitely many of the events $C_x$ occur in both directions, $\Pa$--a.s.
\end{proof}

Assume the environment and cascade are in the event $\mathcal{E}_1$ defiend in the above lemma. Let $(x_k)_{k\in\Z}$ be the locations of coalescence events (with $x_0$ the smallest non-negative $x$ such that $C_x$ occurs). By definition of the $x_k$, for every $k$ and for every $a$ to the left of $[[x_k]]_R$, $H_{[[x_k]]_R}({\bf X}^a)=H_{x_k}({\bf X}^a)<\infty$. Now for $j<k$, it necessarily holds that $x_j$ is to the left of $[[x_k]]_R$, since there can be only one $x_k$ per level. Define $\nu(j,k):=H_{x_k}({\bf X}^{x_j})$. By definition of the walks ${\bf X}^a$, we have for $j<k$, $n\geq0$, 
\begin{equation}\label{eqn:996}
X_{n+\nu(j,k)}^{x_j}=X_n^{x_k}.
\end{equation}
From this one can easily check that the $\nu(j,k)$ are additive; that is, for $j<k<\ell$, we have $\nu(j,\ell)=\nu(j,k)+\nu(k,\ell)$. 
\ifSHOWEXTRA
{\color{blue} We note that for fixed $j$, $\nu(j,k)$ is increasing in $k$, because for $j<k<\ell$, the walk ${\bf X}^j$ must hit $[[x_k]]_R$ before it can hit $[[x_{\ell}]]_R$. 

\begin{align*}
    \nu(j,\ell)&=H_{x_{\ell}}({\bf X}^{x_j})
    \\
    &=\inf\{n\geq0:X_n^{x_j}=x_{\ell}\}
    \\
    &=\nu(j,k)+\inf\{n\geq0:X_{n+\nu(j,k)}^{x_j}=x_{\ell}\}
    \\
    &=\nu(j,k)+\inf\{n\geq0:X_{n}^{x_k}=x_{\ell}\}
    \\
    &=\nu(j,k)+H_{x_k}({\bf X}^{\ell})
    \\
    &=\nu(j,k)+\nu(k,\ell)
\end{align*}
}
\fi
Because all the ${\bf X}^{x_k}$ agree with each other in the sense of \eqref{eqn:996}, we may define a single, bi-infinite walk $\BIX=(\biX_n)_{n\in\Z}$ that agrees with all of the ${\bf X}^{x_k}$. We choose to let $\biX_0=x_0$. For $n\geq0$, let $\biX_n=X_n^{x_0}$. For $n<0$, choose $j<0$ such that $\nu(j,0)>|n|$, and let $X_n=X_{\nu(j,0)-|n|}^{x_j}$. This definition is independent of the choice of $j$, because if $j<k<0$ with $v(k,0)>|n|$, then by \eqref{eqn:996} and the additivity of the $\nu(j,k)$, we have
\begin{equation*}
    X_{\nu(j,0)-|n|}^{x_j}=X_{\nu(j,k)+\nu(k,0)-|n|}^{x_j}=X_{\nu(k,0)-|n|}^{x_k}.
\end{equation*}
%We want to define $\BIX$ to agree with each ${\bf X}^{x_k}$ until it hits the level $[x_{k+1}]$. Let $\tau_k=\inf\{n:X_n^{x_{k-1}}=x_k\}$. These are finite by the definition of the coalescence events. For all $k\in\Z$ and all $0\leq j<\tau_{k+1}-\tau_k$, let $\biX_{\tau_k+j}=X_j^{x_k}$. This defines $\biX_n$ for all $t\in\R$, so we let $\BIX=(\biX_n)_{n\in\Z}$. For any $k$ and any $a\leq x_k$, we note that because ${\bf X}^{a}$ will first hit level $[x_{k+i}]$ at the point $x_{k+i}$ for all $i\geq0$, the walk ${\bf X}^{a}$ agrees with $\BIX$ everywhere to the right of $x_k$.
We may then define $\overline{N}_{x}:=\#\{n\in\Z:\biX_n=x\}$ to be the amount of time the walk $\BIX$ spends at $x$. Thus, $\overline{N}_x=\lim_{a\to-\infty}N_x({\bf X}^a)$.

\begin{lem}\label{prop:Ergodic}
Both of the sequences $({\bf X}^a)_{a\in\Z}$ and $(\overline{N}_x)_{x\in\Z}$ are stationary and ergodic. 
\end{lem}

\begin{proof}
For a given environment, the cascade that defines $\BIX$ may be generated by a (countable) family ${\bf U}=\left(U_n^{a}\right)_{n\in\N,a\in\Z}$ of i.i.d. uniform random variables on $[0,1]$. For such a collection, and an $a\in\Z$, let ${\bf U}^a$ be the projection $\left(U_n^{a}\right)_{n\in\N}$. Given an environment $\omega$, the finite walk from $a$ to level $[[a+R]]_R$ may be generated using the first several $U_n^{a}$. (One of the $U_n^{a}$ is used for each step. Once the walk terminates, the rest of the $U_n^{a}$ are not needed, but one does not know in advance how many will be needed.) Let $\hat{\omega}^x=(\omega^x,{\bf U}^x)$, and $\hat{\omega}=(\hat{\omega}^x)_{x\in\Z}$. Define the left shift $\hat{\theta}$ by $\hat{\theta}(\hat{\omega}):=(\hat{\omega}^{x+1})_{x\in\Z}$.
Then $(\hat{\omega}^x)_{x\in\Z}$ is an i.i.d. sequence. We have ${\bf X}^0={\bf X}^0(\hat{\omega})$ and ${\bf X}^a={\bf X}^0(\hat{\theta}^a\hat{\omega})$. Similarly, $\overline{N}_0=\overline{N}_0(\hat{\omega})$ and $N_x=\overline{N}_0(\hat{\theta}^x\hat{\omega})$. So it suffices to show that ${\bf X}^0$ and $\overline{N}_0$ are measurable.
The measurability of ${\bf X}^0$ is obvious. For $\overline{N}_0$, let $A_{k,\ell,B,r}$ be the event that:
\begin{enumerate}[(a)]
    \item for some $x<0$, a coalescence event $C_{x,k}$ (as defined in the proof of Lemma \ref{prop:coalescences}) occurs with $-B\leq x-kR<x<0$, so that $\BIX$ agrees with ${\bf X}^x$ to the right of $x$;
    \item $N_0^{[-B,B]}({\bf X}^x)\geq \ell$, where $N_0^{[-B,B]}$ is the amount of time the walk spends at $x$ before exiting $[-B,B]$; and
    \item none of the walks from sites $a\in[-B,B]$ uses more than $r$ of the random variables  $U_r^{a}$.
\end{enumerate}
On this event, $\overline{N}_0$ is seen to be at least $\ell$ by looking only within $[-B,B]$ and only at the first $r$ uniform random variables at each site. The event $A_{k,\ell,B,r}$ is measurable, because it is a measurable function of finitely many random variables, and the event $\{\overline{N}_0>\ell\}$ is, up to a null set, simply the union over all $r$, then over all $B$, and then over all $k$ of these events. Thus, $\overline{N}_0$ is measurable.
\end{proof}

We now give the connection between $\overline{N}_0$ and the limiting velocity $v$. 

\begin{lem}\label{lem:LimitingVelocity}
 $v=\frac{1}{\Ea[\overline{N}_0]}$. Consequently, the walk is ballistic if and only if $\Ea[\overline{N}_0]<\infty$.
 \end{lem}

We note that a similar formula for the limiting speed in the ballistic case can be obtained from \cite[Theorem 6.12]{Dolgopyat&Goldsheid2019} for discrete-time RWRE on a strip, but under a stronger ellipticity assumption (and with a less explicit probabilisitic interpretation).

\begin{proof}
By Lemma \ref{prop:Ergodic} and Birkhoff's Ergodic theorem, for any $c\in\Z$ we have
\begin{equation*}\lim_{n\to\infty}\frac1n\sum_{y=c}^n\overline{N}_{y}=\Ea[\overline{N}_0],~\Pa\text{--a.s.}\end{equation*}
Fix $a\in\Z$. For large enough $y$, $N_y({\bf X}^a)=\overline{N}_{y}$. We therefore get
\begin{equation*}\lim_{n\to\infty}\frac1n\sum_{y=c}^nN_y({\bf X}^a)=\Ea[\overline{N}_0],~\Pa\text{--a.s.}\end{equation*}
It follows that
\begin{equation*}\lim_{n\to\infty}\frac1n\sum_{y=c}^nN_y({\bf X})=\Ea[\overline{N}_0],~\Pa^a\text{--a.s.}\end{equation*}
By Lemma \ref{known}, we get $v=\frac{1}{\Ea[\overline{N}_{0}]}$.
\end{proof}

Now we can see that the walk is ballistic if and only if $\Ea[\overline{N}_0]<\infty$. We now compare $\Ea[\overline{N}_0]$ with $\Ea^0[N_0]$. 

\begin{lem}\label{lem:limitfinitecase}
$\Ea[\overline{N}_0]\leq\Ea^0[N_0]$. 
\end{lem}

\begin{proof}
If $\Ea^0[N_0]=\infty$, the inequality is trivial. Assume, therefore, that $\Ea^0[N_0]<\infty$. 

Note that ${\lim_{x\to\infty}N_0({\bf X}^{-x})=\overline{N}_0}$, $\Pa$--a.s. Using Fatou's lemma\footnote{We can get equality in \eqref{eqn:1074} using dominated convergence, but it is a bit more cumbersome and unnecessary.}, we have

\begin{align}
\notag
\Ea[\overline{N}_0] &= \Ea\left[\lim_{x\to\infty}N_0({\bf X}^{-x})\right] 
\\
\label{eqn:1074}
&\leq \lim_{x\to\infty}\Ea\left[N_0({\bf X}^{-x})\right] 
\\
\notag
&= \lim_{x\to\infty}\Ea^{-x}[N_0({\bf X})]. 
\end{align}
But each term $\Ea^{-x}[N_0({\bf X})]=E[E_{\omega}^{-x}[N_0]]$ is less than $\Ea^0[N_0]=E[E_{\omega}^0[N_0]]$, since $E_{\omega}^{-x}[N_0]=P_{\omega}^{-x}(H_0<\infty)E_{\omega}^{0}[N_0]$.
Therefore, we may conclude $\Ea[\overline{N}_0]\leq\Ea^0[N_0]$.
\end{proof}

We are now in a position to prove our main theorem, and in fact two of the three implications are already done. 

\begin{proof}[Proof of Theorem \ref{thm:abstractballisticitycriteria}]
\hspace{1in}

{\em (b) $\Rightarrow$ (c) }
This is an immediate consequence of Lemma \ref{prop:4.8}. 

{\em (c) $\Rightarrow$ (a) }
Assume $\Ea^0[N_0]<\infty$. Then combining Lemmas \ref{lem:LimitingVelocity} and \ref{lem:limitfinitecase} gives $v>0$.

{\em (a) $\Rightarrow$ (b) }
Suppose $\Ea^0[H_{\geq 1}]=\infty$. We will show that $v=0$. We claim
\begin{equation}\label{eqn:509}
    \Ea\left[\min_{1\leq i\leq R}H_{\geq R+1}({\bf X}^i)\right]=\infty.
\end{equation}
By assumptions (\hyperlink{cond:C1}{C1}), (\hyperlink{cond:C2}{C2}), and (\hyperlink{cond:C3}{C3}), there is an $m_0\geq\max(L,R)$ large enough that every interval of length $m_0$ is irreducible with positive $P$-probability.
Let $A$ be the event that
\begin{itemize}
    \item For each $i=1,\ldots, R-1$, the walk ${\bf X}^i$ hits $R$ before leaving $[R-m_0+1, R]$. 
    \item The walk ${\bf X}^R$ first exits $[R-m_0+1,R]$ by hitting $R-m_0$. 
\end{itemize}
Then under $P$, the random variable $P_{\omega}(A)$ is independent of $\omega^{\leq R-m_0}$. Now, on the event $A$, the minimum $\min_{1\leq i\leq R}H_{\geq R+1}({\bf X}^i)$ is attained for $i=R$, since all the other walks take time to get to $R$ and then simply follow ${\bf X}^R$. 
Now on $A$, $H_{\geq R+1}({\bf X}^R)$ is greater than the amount of time it takes for the walk ${\bf X}^R$ to cross back to $[R-m_0+1,\infty)$ after first hitting $R-m_0$. The quenched expectation of this time, conditioned on $A$, is $E_{\omega}^{R-m_0}[H_{\geq R-m_0+1}]$ by the strong Markov property, and this depends only on $\omega^{\leq R-m_0}$. Hence
\begin{align}
\notag
    \Ea\left[\min_{1\leq i\leq R}H_{\geq R+1}({\bf X}^i)\right]
    &\geq E\left[P_{\omega}(A)E_{\omega}[H_{\geq R+1}({\bf X}^R)|A]\right]
    \\
\notag
    &\geq E\left[P_{\omega}(A)E_{\omega}^{R-m_0}[H_{\geq R-m_0+1}({\bf X})]\right]
    \\
\notag
    &= \Pa(A)\Ea^{R-m_0}[H_{\geq R-m_0+1}({\bf X})]
    \\
\notag    
    &=\Pa(A)\Ea^{0}[H_{\geq 1}({\bf X})]
    \\
\notag
    &=\infty.
\end{align}
This proves \eqref{eqn:509}. 
Now for $x\geq1$, 
\begin{equation*}
    H_{\geq xR+1}({\bf X}^0)\geq H_{\geq1}({\bf X}^0)+\sum_{k=1}^x\min_{1\leq i \leq R}H_{\geq kR+1}({\bf X}^{(k-1)R+i}).
\end{equation*}
Dividing by $xR$ and taking limits as $x\to\infty$, we get $\lim_{x\to\infty}\frac{H_{\geq xR+1}({\bf X}^0)}{xR}=\infty$, $\Pa$--a.s. by Birkhoff's ergodic theorem. Hence $\lim_{x\to\infty}\frac{H_{\geq xR+1}({\bf X})}{(x+1)R}=\infty$, $\Pa^0$--a.s. For $n_k=H_{\geq kR+1}({\bf X})$, $X_{n_k}\leq (x+1)R$, and so $\frac{X_{n_k}}{n_k}\leq\frac{(k+1)R}{H_{\geq kR+1}({\bf X})}\to0$ as $k\to\infty$, $\Pa^0$--a.s. Since $\frac{X_{n_k}}{n_k}$ is a subsequence of $\frac{X_n}{n}$, it must $\Pa^0$--a.s. approach $v$, and therefore $v=0$. 
\end{proof}

\bibliographystyle{plain}
\bibliography{default}

\begin{thebibliography}{10}

\bibitem{Bremont2002}
Julien Brémont.
\newblock On some random walks on z in random medium.
\newblock {\em Ann. Probab.}, 30(3):1266--1312, 07 2002.

\bibitem{Bremont2004}
Julien Brémont.
\newblock Random walks in random medium on z and lyapunov spectrum.
\newblock {\em Annales de l'Institut Henri Poincare (B) Probability and
  Statistics}, 40(3):309 -- 336, 2004.

\bibitem{Bremont2009}
Julien Brémont.
\newblock One-dimensional finite range random walk in random medium and
  invariant measure equation.
\newblock {\em Ann. Inst. H. Poincaré Probab. Statist.}, 45(1):70--103, 02
  2009.

\bibitem{Dolgopyat&Goldsheid2019}
Dmitry Dolgopyat and Ilya Goldsheid.
\newblock {Invariant measure for random walks on ergodic environments on a
  strip}.
\newblock {\em The Annals of Probability}, 47(4):2494 – 2528, 2019.

\bibitem{Goldsheid2008}
Ilya Goldsheid.
\newblock Linear and sub-linear growth and the clt for hitting times of a
  random walk in random environment on a strip.
\newblock {\em Probability Theory and Related Fields}, 141:471--511, 07 2008.

\bibitem{Kesten1977}
H.~Kesten.
\newblock A renewal theorem for random walk in a random environment.
\newblock {\em Proc. Sympos. Pure Math.}, 31:67--77, 1977.

\bibitem{Key1984}
Eric~S. Key.
\newblock Recurrence and transience criteria for random walk in a random
  environment.
\newblock {\em Ann. Probab.}, 12(2):529--560, 05 1984.

\bibitem{Roitershtein2008}
Alexander Roitershtein.
\newblock Transient random walks on a strip in a random environment.
\newblock {\em Ann. Probab.}, 36(6):2354--2387, 11 2008.

\bibitem{Slonim2021a}
Daniel~J. Slonim.
\newblock Random walks in dirichlet random environments on $\mathbb{Z}$ with
  bounded jumps.
\newblock Preprint, submitted May 2021. https://arxiv.org/abs/2104.14950.

\bibitem{Sznitman&Zerner1999}
Alain-Sol Sznitman and Martin Zerner.
\newblock {A Law of Large Numbers for Random Walks in Random Environment}.
\newblock {\em The Annals of Probability}, 27(4):1851 – 1869, 1999.

\bibitem{Zeitouni2004}
O.~Zeitouni.
\newblock Random walks in random environment.
\newblock In J.~Picard, editor, {\em Lecture notes in probability theory and
  statistics: École d’été de probabilités de Saint-Flour XXXI-2001},
  volume 1837 of {\em Lect. Notes Math.}, pages 190--313. Springer, 2004.

\bibitem{Zerner2002}
Martin Zerner.
\newblock A non-ballistic law of large numbers for random walks in i.i.d.
  random environment.
\newblock {\em Electron. Commun. Probab.}, 7:191--197, 2002.

\end{thebibliography}

\end{document}